\documentclass[psamsfonts, reqno]{amsart}

\usepackage{amscd,amsmath,amssymb,amsfonts,verbatim}
\usepackage{times}
\usepackage[cmtip, all]{xy}
\usepackage{graphicx}
\usepackage[active]{srcltx}
\usepackage{color}
\usepackage{textcomp}

\definecolor{refkey}{gray}{.5}   
\definecolor{labelkey}{gray}{.5} 
\definecolor{Red}{rgb}{1,0,0}

\newcommand{\pf}{{\bf Proof : }}

\newcommand{\qedwhite}{\hfill \ensuremath{\Box}}

\newtheorem{theo}{Theorem}[subsection]
\newtheorem{prop}[theo]{Proposition}
\newtheorem{lem}[theo]{Lemma}
\newtheorem{cor}[theo]{Corollary}

\theoremstyle{definition}

\newtheorem{rem}[theo]{Remark}
\newtheorem{exam}[theo]{Example}
\newtheorem{defi}[theo]{Definition}

\title{ Absence of torsion in orbit space}
\author{Sampat Sharma}
\newcommand{\Addresses}{{
  \bigskip
  \footnotesize

 \textsc{Sampat Sharma, School of Mathematics, Tata Institute of Fundamental Research,\\ \noindent
           1, Dr. Homi Bhabha Road, Mumbai 400005, INDIA
   } \par\nopagebreak
  \textit{E-mail}: Sampat ~Sharma \texttt{<sampat.iiserm@gmail.com; sampat@math.tifr.res.in>}

  \medskip

  }}

\begin{document}
\maketitle

\vskip0.15in

\subjclass 2020 Mathematics Subject Classification:{13C10, 19B10, 13H99 }

\keywords {Keywords:}~ {Unimodular rows, Suslin matrices}

\begin{abstract}
 In this paper, we prove that if $R$ is a local ring of dimension $d,$ $d\geq 2$ and $\frac{1}{d!}\in R$ then the group $\frac{Um_{d+1}(R[X])}{E_{d+1}(R[X])}$ has no $k$-torsion, provided $kR = R.$ We also prove that if $R$ is a regular ring of dimension 
 $d,$ $d\geq 2$ and $\frac{1}{d!}\in R$ such that $E_{d+1}(R)$ acts transitively on $Um_{d+1}(R)$ then  $E_{d+1}(R[X])$ acts transitively on $Um_{d+1}(R[X]).$
\end{abstract}

\vskip0.50in

\section{Introduction}
Throughout this article, we assume $R$ to be a commutative noetherian ring with $1 \neq 0 $ unless stated otherwise. 
\par In \cite{st}, J. Stienstera, using the ideas of S. Bloch in \cite{bloch} showed that $NK_{1}(R)$ is a $W(R)$-module where $NK_{1}(R) = \mbox{Ker}(K_{1}(R[X])\longrightarrow K_{1}(R)); X = 0$ and $W(R)$ is the ring of big Witt vectors. Consequently, as noted by C. Weibel 
in \cite{weibel}, $SK_{1}(R[X])$ has no $k$-torsion if $kR = R$ and $R$ is a commutative local ring. Note that $SK_{1}(R[X])$ coincides with $NK_{1}(R)$ for a commutative local ring $R$. In \cite{basu}, R. Basu simplified J. Stienstera's approach of big Witt vector and proved that $NK_{1}(R)$ has no $k$-torsion for an associative ring $R$ with $kR =R.$ 

\par In $(${\cite[Theorem 2.2]{fasel}}$),$ J. Fasel proved that for a non-singular affine algebra $R$ of dimension $d\geq 3$ over a perfect field $K$ of cohomological dimension atmost 1, the universal Mennicke symbol $MS_{d+1}(R)$ is uniquely $p$ divisible for $p$ prime to the 
characteristic of $K$. Using results of J. Fasel $(${\cite[Theorem 2.1]{fasel}}$)$ and W. van der Kallen $(${\cite[Theorem 4.1]{vdk2}}$)$, Ravi Rao and Selby Jose in \cite{orthoquot} concluded that for a non-singular affine algebra $R$ of dimension $d\geq 3$ over a perfect field $K$, $\mbox{char}(K)\neq 2,$ of cohomological dimension atmost 1, the orbit space $\frac{Um_{d+1}(R)}{E_{d+1}(R)}$ is uniquely $p$ divisible for $p$ prime to the 
characteristic of $K$.

\par In $(${\cite[Corollary 7.4]{7}}$),$ L.N Vaserstein proves that if $R$ is a local ring of dimension $2,$ 
in which $2$ is invertible then $\frac{Um_{3}(R[X])}{E_{3}(R[X])}\simeq W_{E}(R[X]).$ Here $W_{E}(R[X])$ denotes the elementary symplectic Witt group which has been defined in section 2. Using Karoubi's linearization technique, in \cite{rswan}, Rao-Swan proved that $W_{E}(R[X]) \hookrightarrow SK_{1}(R[X])$  is an injective group homomorphism
 for a local ring $R$ with $2R =R.$ Using Weibel's result in \cite{weibel}, Rao--Swan (for a proof see \cite{orthoquot}), proved that  $\frac{Um_{3}(R[X])}{E_{3}(R[X])}$ has no $2$-torsion 
for a local ring of dimension $2$ with $2R = R.$  

\par In this article, using the theory of Suslin matrices, we generalise the result of Rao--Swan for any $d$-dimensional local ring, $d\geq 2$. We prove that : 
\begin{theo}
\label{introtheo}
 Let $R$ be a local ring of dimension $d$ and let $\frac{1}{d!} \in R$ , then the group
 $\frac{Um_{d+1}(R[X])}{E_{d+1}(R[X])}$ has no $k$-torsion, provided $kR = R.$
\end{theo}
Using Ravi Rao's result $(${\cite[Theorem 2.4]{invent}}$),$ we define a map $\phi : \frac{Um_{d+1}(R[X])}{E_{d+1}(R[X])} \longrightarrow SK_{1}(R[X])$. Weibel's result that $SK_{1}(R[X])$ has no $k$-torsion for a local ring $R$ with $kR = R$, ensures that the 
map $\phi$ is well defined. Next we use a result of Anuradha Garge and Ravi Rao $(${\cite[Corollary 2.4]{gr}}$),$ to establish that $\phi$ is an injective map. Using the fact that $[v]\mapsto [S_{n}(v,w)]$ is a Mennicke symbol, one proves that $\phi$ is a group homomorphism. Here $S_{n}(v,w)$ denotes the Suslin matrix which has been defined in section 2.

\par As an application of the Theorem \ref{introtheo}, we answer a question raised by Ravi Rao and Selby Jose in \cite{orthoquot}. One can define an injective map from orbit space to the quotient group $SO_{2(d+1)}(R[X])/EO_{2(d+1)}(R[X]).$ We prove that 
\begin{cor}
 Let $R$ be a local ring of dimension $d$, $d =2k, k\geq 1$,  $\frac{1}{d!}\in R$. Then the map 
$$\varphi : \frac{Um_{d+1}(R[X])}{E_{d+1}(R[X])}\longrightarrow \frac{SO_{2(d+1)}(R[X])}{EO_{2(d+1)}(R[X])}$$
is injective.
\end{cor}

\par In \cite{4}, A. A. Suslin proved that for a noetherian ring of dimension $d$, $E_{d+2}(R[X])$ acts transitively on $Um_{d+2}(R[X])$ for $d\geq 1.$ As an application of Theorem \ref{introtheo}, we improve this bound upto $d+1$ in certain rings. We prove that : 
\begin{cor} Let $R$ be a regular ring of dimension 
 $d,$ $d\geq 2$ and $\frac{1}{d!}\in R$ such that $E_{d+1}(R)$ acts transitively on $Um_{d+1}(R)$ then  $E_{d+1}(R[X])$ acts transitively on $Um_{d+1}(R[X]).$
\end{cor}

\section{Preliminary Results}
\setcounter{theo}{0}
\par 
A row $v= (a_{0},a_{1},\ldots, a_{r})\in R^{r+1}$ is said to be unimodular if there is a $w = (b_{0},b_{1},\ldots, b_{r})
\in R^{r+1}$ with
$\langle v ,w\rangle = \Sigma_{i = 0}^{r} a_{i}b_{i} = 1$ and $Um_{r+1}(R)$ will denote
the set of unimodular rows (over $R$) of length 
$r+1$.
\par 
The group of elementary matrices is a subgroup of $GL_{r+1}(R)$, denoted by $E_{r+1}(R)$, and is generated by the matrices 
of the form $e_{ij}(\lambda) = I_{r+1} + \lambda E_{ij}$, where $\lambda \in R, ~i\neq j, ~1\leq i,j\leq r+1,~
E_{ij} \in M_{r+1}(R)$ whose $ij^{th}$ entry is $1$ and all other entries are zero. The elementary linear group 
$E_{r+1}(R)$ acts on the rows of length $r+1$ by right multiplication. Moreover, this action takes unimodular rows to 
unimodular
rows : $\frac{Um_{r+1}(R)}{E_{r+1}(R)}$ will denote set of orbits of this action, and we shall denote by $[v]$ the 
equivalence class of a row $v$ under this equivalence relation.

\begin{defi}{\bf{Essential dimension:}}
 Let $R$ be a ring whose maximal spectrum $\mbox{Max}(R)$ is a finite union of subsets $V_{i},$ where each $V_{i},$ 
 when endowed with the (topology induced from the) Zariski topology is a space of Krull dimension $d.$ We shall say 
 $R$ is essentially of dimension $d$ in such a case.
\end{defi}
\par 
In $(${\cite[Theorem 3.6]{vdk1}}$),$ W. van der Kallen derives an abelian
 group structure on $\frac{Um_{d+1}(R)}{E_{d+1}(R)}$ when 
$R$ is of essential dimension $d,$ for all $d\geq 2.$ We will denote the group operation in this group by $\ast.$

\subsection{Abelian group structure on $\frac{Um_{n}(R)}{E_{n}(R)}$ for $\mbox{sdim}(R)\leq 2n-4$}
\begin{defi}{{Stable range condition}}~{{$sr_{n}(R):$}} We shall say stable range condition $sr_{n}(R)$ holds for an associative ring $R$ if 
 for any $(a_{1},a_{2},\ldots,a_{n+1})\in Um_{n+1}(R)$ there exists $c_{i}\in R$ such that 
 $$(a_{1}+ c_{1}a_{n+1}, a_{2}+c_{2}a_{n+1},\ldots,a_{n}+c_{n}a_{n+1})\in Um_{n}(R).$$
 
\end{defi}
\begin{defi}{Stable range $sr(R),$ Stable dimension $\mbox{sdim}(R):$} The stable range
of an associative ring $R$ is defined to be the least integer $n$ such that $sr_{n}(R)$ holds. Stable range of $R$ is denoted by $sr(R).$ The stable dimension of $R$ is defined by
$\mbox{sdim}(R) = sr(R) - 1.$
 
\end{defi}

\begin{defi}{Universal weak Mennicke symbol $\mbox{WMS}_{n}(R), n\geq 2:$} We define the universal weak Mennicke symbol on 
$\frac{Um_{n}(R)}{E_{n}(R)}$ by a set map $\textit{wms} : \frac{Um_{n}(R)}{E_{n}(R)} \longrightarrow \mbox{WMS}_{n}(R),$ 
$[v]\longmapsto
\textit{wms}(v)$
 to a group $\mbox{WMS}_{n}(R).$ The group $\mbox{WMS}_{n}(R)$ is the free abelian group generated by $\textit{wms}(v),
 v\in Um_{n}(R)$ modulo the 
 following relations 
 \begin{itemize}
  \item $\textit{wms}(v) = \textit{wms}(v\varepsilon)~\mbox{if}~\varepsilon \in E_{n}(R).$
  \item If $(q, v_{2}, \ldots, v_{n}), (1+q, v_{2}, \ldots, v_{n}) \in Um_{n}(R)$ and $r(1+q) \equiv q~\mbox{mod} 
  ~(v_{2}, \ldots
  , v_{n}), $ then 
  $$\textit{wms}(q, v_{2}, \ldots, v_{n}) = \textit{wms}(r, v_{2}, \ldots, v_{n})\textit{wms}(1+q, v_{2}, \dots, v_{n}),~~~~\mbox{for}~r\in R.$$
 \end{itemize}

\end{defi}
We recall $(${\cite[Theorem 4.1]{vdk2}}$).$
\begin{theo}(W. van der Kallen) Let $R$ be a ring of stable dimension $d,$ $d\leq 2n-4$ and $n\geq 3.$ Then the universal 
weak Mennicke symbol $\textit{wms} : \frac{Um_{n}(R)}{E_{n}(R)} \longrightarrow \mbox{WMS}_{n}(R)$ is bijective.
 
\end{theo}

\par In view of above theorem $\frac{Um_{n}(R)}{E_{n}(R)}$ gets a group structure for $d\leq 2n-4, 
n\geq 3$. We will denote group operation in 
$\frac{Um_{n}(R)}{E_{n}(R)}$ by $\ast.$ In $(${\cite[Lemma 3.5]{vdk2}}$),$ W. van der Kallen proves various revealing
formulae in this group. Let us recall some of these; 
\begin{itemize}
 \item $[(x_{1}, v_{2}, \ldots, v_{n})]\ast [(v_{1}, v_{2}, \ldots, v_{n})] = [(v_{1}(x_{1}+ w_{1})-1, (x_{1}+w_{1})v_{2}, 
 \ldots, v_{n})],$ where $w_{1}$ is such that $\Sigma_{i=1}^{n}v_{i}w_{i} = 1,$ for some $w_{i}\in R,$ $1\leq i\leq n.$ In 
 particular $[(v_{1}, v_{2}, \ldots, v_{n})]^{-1} = [(-w_{1}, v_{2}, \ldots, v_{n})].$
 \item $[(x_{1}, v_{2}, \ldots, v_{n})]\ast [(v_{1}, v_{2}, \ldots, v_{n})]^{-1} = [(1-w_{1}(x_{1}-v_{1}), (x_{1}-v_{1})v_{2}, 
 \ldots, v_{n})].$
 \item $[(x_{1}, v_{2}, \ldots, v_{n})]\ast [(v_{1}^{2}, v_{2}, \ldots, v_{n})] = [(x_{1}v_{1}^{2}, v_{2}, \ldots, v_{n})].$
\end{itemize}

\subsection{Mennicke--Newman Lemma}
We note Mennicke--Newman lemma proved by W. van der Kallen in $(${\cite[Lemma 3.2]{vdk3}}$)$: 
\begin{lem}
\label{mennickenewman}
 Let $R$ be a commutative ring with $\mbox{sdim}(R) \leq 2n-3.$ Let $v, w\in Um_{n}(R).$ Then there exists 
 $\varepsilon, \delta \in E_{n}(R)$ and $x, y, a_{i}\in R$ such that for $i =2, \ldots, n$
 $$v\varepsilon = (x, a_{2}, \ldots, a_{n}), w\delta = (y, a_{2}, \ldots, a_{n}),~ x+y = 1.$$
\end{lem}

 \subsection{L.N Vaserstein's power operation $\chi_{k}$ on $\frac{Um_{r}(R)}{E_{r}(R)},$ for $k\in \mathbb{Z}, r\geq 3$}

 In \cite{orbit}, L.N. Vaserstein has shown that taking $k^{th}$ power of a co-ordinate is a well defined operation $\chi_{k}$ on 
$\frac{Um_{r}(R)}{E_{r}(R)}, r\geq 3,$ for any commutative ring $R.$ Let $\overset{E}\sim$ denote equivalence under the 
elementary group $E_{r}(A).$ Thus, 
\begin{itemize}
 \item $(v_{1},\ldots,v_{i}^{k},\ldots, v_{r}) \overset{E}\sim (v_{1},\ldots,v_{j}^{k},\ldots, v_{r}),$ for
 all $1\leq i, j \leq
  r, k\geq 0.$
  \item $(v_{1},\ldots,v_{i-1},w_{i},v_{i+1},\ldots, v_{r}) \overset{E}\sim (v_{1},\ldots,v_{j-1},w_{j},v_{j+1},
  \ldots, v_{r}),$
   for all $1\leq i, j \leq r,$ where $\Sigma_{i=1}^{r}v_{i}w_{i} = 1,$ for some $w_{i}\in R.$
\end{itemize}

\subsection{On the equality of $\chi_{k}([v])$ with $[v]^{k}, k\in \mathbb{N}$}
\setcounter{theo}{0}
W. van der Kallen in $(${\cite[Section 4]{vdk2}}$),$ gives an example to illustrate that, in general,
$(v_{1}^{k}, v_{2}, \ldots
, v_{n}) \notin [v]^{k},$ where $v = (v_{1}, \ldots, v_{n}).$  A cause for this anomaly was 
pointed out by Ravi A. Rao in {\cite{invent}} $-$ Antipodal unimodular 
rows may not coincide up to an elementary action. Following $(${\cite[Lemma 1.3.1]{invent}}$),$ we prove the following : 
\begin{lem}
\label{antipod}
Let $v = (v_{1}, v_{2}, \ldots, v_{n})\in Um_{n}(R),$ where $sdim(R) = d,$ $d\leq 2n-4.$ Let us assume that  
$v \overset{E}\sim (-v_{1}, v_{2}, \ldots, v_{n})$ (equivalently, $\chi_{-1}([v]) = [v]^{-1}$), then 
$(v_{1}^{k}, v_{2}, \ldots
, v_{n}) \in [v]^{k}$ for all $k\in \mathbb{N}.$
\end{lem}
${\pf}$ Let $\Sigma_{i=1}^{n}v_{i}w_{i} = 1,$ for some $w_{i}\in R,$ $1\leq i\leq n.$ Note that we can then write 
$v_{1}^{2}w_{1}^{2} + \Sigma_{i=2}^{n}v_{i}w_{i}' = 1$ for some $w_{i}'\in R,$ and $2\leq i\leq n.$ Now 
$$[(w_{1}, v_{2}, \ldots, v_{n})] = \chi_{-1}([v]) = \chi_{-1}([(-v_{1}, \ldots, v_{n})]) = [(-w_{1}, v_{2}, \ldots, v_{n})].$$
Hence, 
\begin{align*}
[v] = [(-w_{1}, v_{2}, \ldots, v_{n})]^{-1} & = [(v_{1}^{2}w_{1}, v_{2}, \ldots, v_{n})]\\
& = [(w_{1}, v_{2}, \ldots, v_{n})]\ast [(v_{1}^{2}, v_{2}, \ldots, v_{n})]\\
& = [(-w_{1}, v_{2}, \ldots, v_{n})] \ast [(v_{1}^{2}, v_{2}, \ldots, v_{n})]\\
& = [v]^{-1} \ast [(v_{1}^{2}, v_{2}, \ldots, v_{n})].
\end{align*}
Therefore $[v]^{2} = [(v_{1}^{2}, v_{2}, \ldots, v_{n})].$ It is now easy to check, via the group identities in 
$\frac{Um_{n}(R)}{E_{n}(R)}$ that\\ $[(v_{1}^{k}, v_{2}, \ldots
, v_{n})] = [v]^{k},$ for all $k > 2.$
$~~~~~~~~~~~~~~~~~~~~~~~~~~~~~~~~~~~~~~~~~~~~~~~~~~~~~~~~~~~~~~~~~~~~~~~~~~~~~~~~~~~~~~~~~~~~~~~~~~~~~~~~~~~~~~\qedwhite$

\subsection{The elementary symplectic Witt group $W_{E}(R)$} 
If $\alpha \in M_{r}(R), \beta \in M_{s}(R)$ are matrices 
then $\alpha \perp \beta$ denotes the matrix $\begin{bmatrix}
                 \alpha & 0\\
                 0 & \beta\\
                \end{bmatrix} \in M_{r+s}(R)$. $\psi_{1}$ will denote $\begin{bmatrix}
                 0& 1\\
                -1 & 0\\
                \end{bmatrix} \in E_{2}(\mathbb{Z})$, and $\psi_{r}$ is inductively defined by $\psi_{r} = \psi_{r-1}\perp 
                \psi_{1} \in E_{2r}(\mathbb{Z})$, for $r\geq 2$.
                \par 
    A skew-symmetric matrix whose diagonal elements are zero is called an alternating matrix. If $\phi \in M_{2r}(R)$ is 
    alternating then $\mbox{det}(\phi) = (\mbox{pf}(\phi))^{2},$ where $\mbox{pf}$ is a 
    polynomial (called the Pfaffian) in the matrix elements 
    with coefficients $\pm 1$. Note that we need to fix a sign in the choice of $\mbox{pf},$ and hence we insist 
    $\mbox{pf}(\psi_{r}) = 1$ 
    for all $r$. For any $\alpha \in M_{2r}(R)$ and any alternating matrix $\phi \in M_{2r}(R)$ we have 
    $\mbox{pf}(\alpha^{t}\phi
    \alpha) = \mbox{pf}(\phi)\mbox{det}(\alpha)$. For alternating matrices $\phi, \psi$ it is easy to check that 
    $\mbox{pf}(\phi \perp \psi) 
    = (\mbox{pf}(\phi))(\mbox{pf}(\psi))$.
    \par Two matrices $\alpha \in M_{2r}(R), \beta \in M_{2s}(R)$ are said to be equivalent (w.r.t. $E(R)$) if there exists  
    a matrix $\varepsilon \in SL_{2(r+s+l)}(R) \bigcap E(R)$,  such that $\alpha \perp \psi_{s+l} = \varepsilon^{t}
    (\beta \perp \psi_{r+l})\varepsilon,$ for some $l$. Denote this by $\alpha \overset{E}\sim \beta$. Thus $\overset{E}\sim $ is an 
    equivalence relation; denote by $[\alpha]$ the orbit of $\alpha$ under this relation. 
    
    \par It is easy to see $(${\cite[p. 945]{7}}$)$ that $\perp$ induces the structure of an abelian group on the set of all 
    equivalence classes of alternating matrices with Pfaffian $1.$ This group is called elementary symplectic Witt group 
    and is denoted by $W_{E}(R)$.

   \subsection{The Suslin matrices}
   First recall the Suslin matrix $S_{r}(v,w)$. These were defined by Suslin in 
$(${\cite[Section 5] {sus2}}$)$. We recall his inductive process : Let $v = (a_{0}, v_{1}), w = (b_{0}, w_{1})$, where 
$v_{1}, w_{1} \in M_{1,r}(R)$. Set $S_{0}(v,w) = a_{0}$ and set 
$$S_{r}(v,w) = \begin{bmatrix}
                 a_{0}I_{2^{r-1}} & S_{r-1}(v_{1},w_{1})\\
                 -S_{r-1}(w_{1}, v_{1})^{T} & b_{0}I_{2^{r-1}}\\
                \end{bmatrix}.$$

\par 
The process is reversible and given a Suslin matrix $S_{r}(v,w)$ one can recover the associated rows $v,w$, i.e. the pair 
$(v,w)$. 

\par 
In \cite{suskcohomo}, Suslin proves that if $\langle v ,w\rangle = v\cdot w^{T} = 1$, then by row and 
column operations one can reduce 
$S_{r}(v,w)$ to a matrix $\beta_{r}(v,w)\in SL_{r+1}(R),$ whose first row is $(a_{0},a_{1},a_{2}^{2}\cdots,a_{r}^{r})$. 
This in particular proves that rows of such type can be completed to an invertible matrix of determinant one. We call any such
$\beta_{r}(v,w)$ to be a compressed Suslin matrix.

\subsection{Relative elementary groups}
\begin{defi}
     Let $I$ be an ideal of a ring $R.$ A unimodular row $v\in Um_{n}(R)$ which is congruent to $e_{1}$ modulo $I$ is 
     called unimodular relative to ideal $I.$ Set of unimodular rows relative to ideal $I$ will be denoted by $Um_{n}(R,I).$
    \end{defi}

Let $I$ be an ideal of a ring $R$, we shall denote by $GL_{n}(R,I)$ the kernel of the canonical mapping 
$GL_{n}(R)\longrightarrow GL_{n}\left({R}/{I}\right).$ Let $SL_{n}(R,I)$ denotes the subgroup of
$GL_{n}(R,I)$ consisting of elements of determinant $1$.

\begin{defi}${\bf{The~ relative~ groups~ E_{n}(I),~ E_{n}(R,I):}}$
 Let $I$ be an ideal of $R$. The elementary group $E_{n}(I)$ is the subgroup of $E_{n}(R)$ generated
 as a group by the elements
 $e_{ij}(x),~x\in I,~1\leq i\neq j\leq n.$\\
  The relative elementary group $E_{n}(R,I)$ is the normal closure of $E_{n}(I)$ in $E_{n}(R)$.
  
\end{defi} 

  \begin{defi}${\bf {Excision ~ring:}}$
 Let $R$ be a ring and $I$ be an ideal of $R$. The excision ring $R\oplus I$, has coordinate wise addition and multiplication
 are given as follows:
 $$(r,i).(s, j) = (rs, rj+si+ij),~\mbox{where}~r,s \in R ~\mbox{and}~i,j\in I.$$
 \par
 The multiplicative identity of this group is $(1,0)$ and the additive identity is $(0,0)$.
\end{defi}

\begin{lem}
\label{anjan}
 $(${\cite[Lemma 4.3]{ggr}}$)$
 Let $(R, \mathfrak{m})$ be a local ring. Then the excision ring $R\oplus I$ with respect to a proper ideal $I\subsetneq R$ is 
 also a local ring with maximal ideal $\mathfrak{m}\oplus I$.
\end{lem}         
 
 \begin{defi}
  For two commutative rings $B$ and $D,$ we say a ring homomorphism $\phi: B\twoheadrightarrow D$ has a section if there exists a ring homomorphism 
  $\gamma : D \hookrightarrow B$ so that $\phi \circ \gamma $ is the identity on $D.$ We shall also say $D$ is a retract 
  of $B.$
 \end{defi}
 
 The following lemma is an easy consequence of $(${\cite[Lemma 4.3, Chapter 3]{mg}}$):$

 \begin{lem}
  $($A. Suslin$)$ 
  \label{k2}
  Let $B,D$ be rings and $\pi: B \twoheadrightarrow D$ be a section. If $J = ker(\pi),$ then 
  $E_{n}(B,J) = E_{n}(B) \cap SL_{n}(B,J), ~n\geq 3.$
 \end{lem}

 \begin{rem}
\label{remark}
   Let $R$ be a ring and $I$ be an ideal of $R$. There is a natural homomorphism $\omega : R\oplus I \rightarrow R$ 
   given by $(x,i) \mapsto x +i \in R.$ Clearly $\omega$ has a section. Let 
   $v = (1+i_{1}, i_{2}, \ldots, i_{n}) \in Um_{n}(R,I)$ where $i_{j}$'s are in $I.$ Then we shall call 
  $\tilde{v} = ((1,i_{1}), (0,i_{2}), \ldots, (0,i_{n}))\in Um_{n}(R\oplus I, 0\oplus I)$ to be a lift of $v.$ Note 
  that $\omega$ sends $\tilde{v}$ to $v.$
 \end{rem}

\section{Absence of torsion in $\frac{Um_{d+1}(R[X])}{E_{d+1}(R[X])}$}
\setcounter{theo}{0}
In $(${\cite[Corollary 7.4]{7}}$),$ L.N Vaserstein proves that if $R$ is a local ring of dimension $2,$ 
in which $2$ is invertible then $\frac{Um_{3}(R[X])}{E_{3}(R[X])}\simeq W_{E}(R[X]).$ Since, $W_{E}(R[X]) \hookrightarrow SK_{1}(R[X])$  (See  Corollary \ref{rao-swan}) for a local ring $R$ with $2R =R,$ one gets  $\frac{Um_{3}(R[X])}{E_{3}(R[X])}$ has no $2$-torsion 
for a local ring of dimension $2.$ In this section we generalise this 
result for any $d$-dimension local ring, $d\geq 2$, via a different approach.
 \par We first collect some known results which will be used in this section.
           
           \begin{prop}
           $(${\cite[Section 3]{weibel}}$)$
           \label{wiebel}
            Let $R$ be a local ring. Assume that $kR = R$. Then $SK_{1}(R[X])$ has no $k$-torsion.
           \end{prop}
\par
Next we recall an observation of Rao-Swan in \cite{rswan}, about a result of Karoubi, and its consequence, proofs can be found in 
 $(${\cite[Lemma 4.2, Corollary 4.3]{gr}}$)$.
 
 \begin{lem}
  {\bf{(Karoubi)} } Let $\phi + NX$, with $\phi \in GL_{2k}(R), N\in M_{2k}(R)$, be 
  a linear invertible alternating polynomial matrix over $R$. If $2R = R$, then $\phi + NX = W(X)^{t}(\phi^{-1})^{t}W(X)$, 
  where $W(X) = \phi(\textit{Id}+\phi^{-1}NX)^{\frac{1}{2}}$.
 \end{lem}
\begin{cor}
\label{rao-swan}
 {\bf{(Rao-Swan)} } Let $R$ be a local ring with $1$, in which $2R = R$. Then the natural map $W_{E}(R[X])\longrightarrow 
 SK_{1}(R[X])$ is an injective group homomorphism.
\end{cor}
${\pf}$ For a proof see $(${\cite[Corollary 4.3]{gr}}$).$
$~~~~~~~~~~~~~~~~~~~~~~~~~~~~~~~~~~~~~~~~~~~~~~~~~~~~~~~~~~~~~~~~~~~~~~~~~~~~~~~~~~~~~~~~~~~~~~~~~~~~~~~~~
           ~~~~~\qedwhite$

\begin{cor}
 Let $R$ be a local ring and $kR = R$. Then $W_{E}(R[X])$ does not have any $k$-torsion.
\end{cor}
${\pf}$ It follows from Proposition \ref{wiebel} and Corollary \ref{rao-swan}.
$~~~~~~~~~~~~~~~~~~~~~~~~~~~~~~~~~~~~~~~~~~~~~~~~~~~~~~~~~~~~~~~~~~~~~~~~~~~~~~~~~~~~~~~~~~~~~~~~~~~~~~~~~
           ~~~~~\qedwhite$

\begin{theo}
\label{factcomp}
$(${\cite[Theorem 2.4]{invent}}$)$
 Let $R$ be a local ring of dimension $d\geq 2$. Let $v\in Um_{d+1}(R[X])$. If $\frac{1}{d!}\in R$ then 
 $v \overset{E}{\sim} (w_{0}^{d!},w_{1},\ldots,w_{d})$ for some $(w_{0}, \ldots, w_{d})\in Um_{d+1}(R[X])$.
\end{theo} 
\begin{rem}
 \label{generalfactcomp}
 In general, if $R$ is a local ring of dimension $d,$ $d\geq 2,$ $2mR = R$ and $v\in Um_{d+1}(R[X]).$ Then 
 one can similarly prove that 
  $$v\overset{E_{d+1}(R[X])}{\sim} (u_{1}^{2m}, u_{2},\ldots, u_{d+1}) $$ for some 
            $(u_{1},\ldots, u_{d+1})\in Um_{d+1}(R[X]).$
 \end{rem}

\begin{theo}
\label{torsion1}
 Let $R$ be a commutative ring with $1,$ $\mbox{sdim}(R) \leq 2n-2,$ $\frac{1}{n!}\in R$ and $n\geq 2.$ Assume that 
 \begin{enumerate}
  \item For all $v\in Um_{n+1}(R)$, $v\overset{E_{n+1}(R)}{\sim} (w_{0}^{n!}, w_{1}, \ldots, w_{n}).$
  \item If $\beta_{n}(v,w)\in SL_{n+1}(R) \cap E_{m}(R), m\geq n+2$ then $[e_{1}\beta_{n}(v,w)] = [e_{1}].$
  \item $SK_{1}(R)$ has no $k$-torsion, if $kR = R.$
 \end{enumerate}
Then, for some $w$ with $\langle e_{1}\sigma, w\rangle = 1, \sigma \in SL_{n+1}(R),$ the map   
\begin{align*}
             \varphi :  &\frac{Um_{n+1}(R)}{E_{n+1}(R)} \longrightarrow SK_{1}(R)\\
             &[v] = [e_{1}\sigma]^{n!}\longmapsto [S_{n}(e_{1}\sigma,w)]
            \end{align*}
is a well-defined injective group homomorphism. Hence  the group $\frac{Um_{n+1}(R)}{E_{n+1}(R)}$ has no $k$-torsion.
\end{theo}
${\pf}$ {\bf{Well-defined:}} In view of the hypothesis $(1)$, for any $v\in Um_{n+1}(R)$, 
 there exists $v'\in Um_{n+1}(R)$, such that $[v] = \chi_{n!}([v']).$ Since every $v'\in Um_{n+1}(R)$ 
is completable, so by Lemma \ref{antipod}, $[v] = [e_{1}\sigma]^{n!},$ for 
some $\sigma \in SL_{n+1}(R).$ 
\par If $v', v'' \in Um_{n+1}(R)$ be such that $[v] = [v']^{n!} = [v'']^{n!}.$ By  
hypothesis $(1)$, there exists $\tau, \tau' \in SL_{n+1}(R)$ such that $[v] = [e_{1}\tau]^{n!} = [e_{1}\tau']^{n!}.$ 
By 
$(${\cite[Lemma 3.2]{key}}$)$, $[S_{n}(e_{1}\tau^{n!},w')] = [S_{n}(e_{1}\tau'^{n!},w'')].$ Since $[v]\longmapsto [S_{n}(v,w)]$ 
is a Mennicke symbol and using Lemma \ref{antipod}, we have
$[S_{n}(e_{1}\tau,w')]^{n!} = [S_{n}(e_{1}\tau',w'')]^{n!}.$ Since $SK_{1}(R)$ does not have any $(n!)$-torsion, 
we have $[S_{n}(e_{1}\tau,w')] = [S_{n}(e_{1}\tau',w'')].$
Therefore the map 
\begin{align*}
             \varphi :  &\frac{Um_{n+1}(R)}{E_{n+1}(R)} \longrightarrow SK_{1}(R)\\
             &[v] = [e_{1}\sigma]^{n!}\longmapsto [S_{n}(e_{1}\sigma,w)]
            \end{align*}
for some $w\in Um_{n+1}(R)$ with $(e_{1}\sigma)\cdot  w^{t} = 1,$ is well-defined.
\par
{\bf{Injectivity:}} Let $v = [e_{1}\sigma]^{n!}$ be such that $[S_{n}(e_{1}\sigma,w)] = 1$ in $SK_{1}(R).$ By 
$(${\cite[Theorem 3.2]{V}}$)$, $ \beta_{n}(e_{1}\sigma, w) \in SL_{n+1}(R) \cap E(R).$ By 
hypothesis $(2),$ 
$[v] = [e_{1}\sigma]^{n!} =  [e_{1}].$ Thus $\phi$ is injective. 
\par 
{\bf{Homomorphism:}} Let $v, w \in Um_{n+1}(R)$ be such that $[v] = [e_{1}\sigma]^{n!}$ and $[w] = [e_{1}\tau]^{n!}$ for 
some matrices $\sigma , \tau \in SL_{n+1}(R).$ In view of Mennicke-Newman lemma (Lemma \ref{mennickenewman}), 
we may assume that 
$[e_{1}\sigma] = [a, x_{1}, x_{2}, \cdots,x_{n}]$ and $[e_{1}\tau] = [b, x_{1}, x_{2}, \cdots,x_{n}].$ Now, in view of 
Lemma \ref{antipod}, we have $[e_{1}\sigma]^{n!} = [a^{n!}, x_{1}, x_{2}, \cdots,x_{n}]$ and 
$[e_{1}\tau]^{n!} = [b^{n!}, x_{1}, x_{2}, \cdots,x_{n}].$ Thus by $(${\cite[Theorem 3.16(iii)]{vdk1}}$),$ we have 
$[v]*[w] = [(ab)^{n!}, x_{1}, x_{2}, \cdots,x_{n}].$
\par Now, by hypothesis $(1),$ there exists $\tau' \in SL_{n+1}(R)$ 
be such that $[e_{1}\tau'] = [ab, x_{1}, x_{2}, \cdots,x_{n}].$ Thus, 
\begin{align*}
             \varphi([v]*[w]) &= [S_{n}(e_{1}\tau' , w')]\\
             & = [S_{n}((ab,x_{1},\ldots,x_{n}) , w')].
            \end{align*}
            \par Since $v \longmapsto S_{n}(v,w)$ is a Mennicke symbol, we have 
            \begin{align*}
             \varphi([v]*[w]) & = [S_{n}((ab,x_{1},\ldots,x_{n}) , w')]\\
             & = [S_{n}((a,x_{1},\ldots,x_{n}) , w_{1})][S_{n}((b,x_{1},\ldots,x_{n}) , w_{2})]\\
             & = [S_{n}(e_{1}\sigma , w_{1})][S_{n}(e_{1}\tau , w_{2})]\\
             & =  \varphi([v])\varphi([w]).
            \end{align*}
It shows that there is an injective group homomorphism $$\frac{Um_{n+1}(R)}{E_{n+1}(R)}\hookrightarrow SK_{1}(R).$$ Since $SK_{1}(R)$ 
does not have any 
$k$-torsion, $\frac{Um_{n+1}(R)}{E_{n+1}(R)}$ has no $k$-torsion.
$~~~~~~~~~~~~~~~~~~~~~~~~~~~~~~~~~~~~~~~~~~~~~~~~~~~~~~~~~~~~~~~~~~~~~~~~~~~~~~~~~~~~~~~~~~~~~~~~~~~~~~~~~
           ~~~~~\qedwhite$

\begin{cor}
\label{absence}
 Let $R$ be a local ring of dimension $d, d\geq 2$ and let $\frac{1}{d!} \in R$ , then 
 $\frac{Um_{d+1}(R[X])}{E_{d+1}(R[X])}$ has no $k$-torsion, provided $kR = R.$
\end{cor}
${\pf}$ It follows from Theorem \ref{torsion1} as all the hypothesis are satisfied in view of Theorem \ref{factcomp}, 
$(${\cite[Corollary 2.2]{gr}}$),$ and Proposition \ref{wiebel} respectively.
$~~~~~~~~~~~~~~~~~~~~~~~~~~~~~~~~~~~~~~~~~~~~~~~~~~~~~~~~~~~~~~~~~~~~~~~~~~~~~~~~~~~~~~~~~~~~~~~~~~~~~~~~~
           ~~~~~\qedwhite$

\begin{lem}
   \label{relativeabsence}
   Let $R$ be a local ring of dimension $d\geq 2,$  $\frac{1}{d!}\in R$ and $I$ be a proper ideal in $R.$
   Let $v\in Um_{d+1}(R[X], I[X])$. Then $\frac{Um_{d+1}(R[X], ~I[X])}{E_{d+1}(R[X],~ I[X])}$ has no $k$-torsion, 
    provided $kR = R.$
  \end{lem}
  ${\pf}$ Let us assume that $[v]^{k} \overset{E_{d+1}(R[X], I[X])}{\sim} e_{1}.$ 
  Let $\tilde{v} \in Um_{d+1}((R\oplus I)[X], (0\oplus I)[X]).$ Here $\tilde{v}$ is defined according to Remark \ref{remark}. In view of Lemma \ref{anjan} and 
  $($ {\cite[Lemma 3.19]{vdk1}}$),$ $R\oplus I$ is a local ring of dimension $d.$ By Corollary \ref{absence}, 
  the group $\frac{Um_{d+1}((R\oplus I)[X])}{E_{d+1}((R\oplus I)[X])}$ has no $k$-torsion. 
  Thus there exists $\varepsilon\in E_{d+1}((R\oplus I)[X])$ such that $\tilde{v}\varepsilon = e_{1}.$ 
  Going modulo $0\oplus I,$ we have $e_{1}\overset{-}\varepsilon = e_{1}, \overset{-}\varepsilon \in 
  SL_{d+1}(R[X]).$ Now replacing $\varepsilon$ by $\varepsilon(\overset{-}\varepsilon)^{-1}$ and
  using 
  Lemma \ref{k2}, we may assume that $\varepsilon \in E_{d+1}((R\oplus I)[X], (0\oplus I)[X])$ satisfying 
  $\tilde{v}\varepsilon = e_{1}.$ Now applying $\omega$ to last equation we get $v\varepsilon_{1} = e_{1}$ for some 
  $\varepsilon_{1} \in E_{d+1}(R[X], I[X])$ which proves that $\frac{Um_{d+1}(R[X],~ I[X])}{E_{d+1}(R[X], ~I[X])}$ has no 
   $k$-torsion.
     $~~~~~~~~~~~~~~~~~~~~~~~~~~~~~~~~~~~~~~~~~~~~~~~~~~~~~~~~~~~~~~~~~~~~~~~~~~~~~~~~~~~~~~~~~~~~~~~~~~~~~~~~~
           ~~~~~\qedwhite$  

\begin{cor}
 Let $R$ be a local ring of dimension $d$, $d =2m, m\geq 1$,  $\frac{1}{d!}\in R$. Then the map 
$$\varphi : \frac{Um_{d+1}(R[X])}{E_{d+1}(R[X])}\longrightarrow \frac{SO_{2(d+1)}(R[X])}{EO_{2(d+1)}(R[X])}$$
is injective.
\end{cor}

 ${\pf}$ Let $[v] \in  \frac{Um_{d+1}(R[X])}{E_{d+1}(R[X])}$ such that $\varphi([v]) \in EO_{2(d+1)}(R[X]).$ In view of $($ {\cite[Theorem 3.8]{orthoquot}}$),$ $\chi_{2}[v] = e_{1}.$ Now,  we use Corollary \ref{absence} to conculde that $[v] = e_{1},$ i.e. $\varphi$ is injective.
 $~~~~~~~~~~~~~~~~~~~~~~~~~~~~~~~~~~~~~~~~~~~~~~~~~~~~~~~~~~~~~~~~~~~~~~~~~~~~~~~~~~~~~~~~~~~~~~~~~~~~~~~~~
           ~~~~~\qedwhite$

\section{An application to completion of unimodular rows}
In $($ {\cite[Theorem 2.6]{4}}$),$ A. A. Suslin proved the following result: 
\begin{theo}$($A.A. Suslin$)$ 
Let $R$ be a noetherian ring of dimension $d$, then $E_{r}(R[X])$ acts transitively on $Um_{r}(R[X])$ for $r\geq \mbox{Max}(3, d+2)$.
\end{theo}

As an application of Theorem \ref{torsion1}, we improve the bound to $d+1$ over regular local rings. We prove : 
\begin{cor}\label{regular}  Let $R$ be a regular local ring of dimension $d, d\geq 2$ and let $\frac{1}{d!} \in R$ , then $E_{r}(R[X])$ acts transitively on $Um_{r}(R[X])$ for $r\geq \mbox{Max}(3, d+1)$.
\end{cor}
${\pf}$ If $r\geq d+2$ then it follows from Suslin's theorem. Let $r = d+1.$ In view of theorem \ref{torsion1}, there is an injective group homomorphism  
$$\varphi : \frac{Um_{d+1}(R[X])}{E_{d+1}(R[X])}\longrightarrow SK_{1}(R[X]).$$
Since $R$ is a regular local ring $SK_{1}(R[X])=0.$ Thus $E_{d+1}(R[X])$ acts transitively on $Um_{d+1}(R[X]).$ 
$~~~~~~~~~~~~~~~~~~~~~~~~~~~~~~~~~~~~~~~~~~~~~~~~~~~~~~~~~~~~~~~~~~~~~~~~~~~~~~~~~~~~~~~~~~~~~~~~~~~~~~~~~
           ~~~~~\qedwhite$

\begin{cor}\label{d+1} Let $R$ be a regular ring of dimension $d, d\geq 2$, $\frac{1}{d!} \in R$  such that $E_{d+1}(R)$ acts transitively on $Um_{d+1}(R)$, then $E_{r}(R[X])$ acts transitively on $Um_{r}(R[X])$ for $r\geq \mbox{Max}(3, d+1)$. 
\end{cor}
${\pf}$ If $r\geq d+2$ then it follows from Suslin's theorem.  Let $r = d+1.$ For any maximal ideal $\mathfrak{m}$ of $R$, in view of corollary \ref{regular}, $E_{d+1}(R_{\mathfrak{m}}[X])$ acts transitively on $Um_{d+1}(R_{\mathfrak{m}}[X]).$ By Sulsin's local-global principle, we gets $v(X)\underset{E_{d+1}(R[X])}{\sim} v(0)$ for every $v(X)\in Um_{d+1}(R[X]).$ Now result follows from our assumption.
$~~~~~~~~~~~~~~~~~~~~~~~~~~~~~~~~~~~~~~~~~~~~~~~~~~~~~~~~~~~~~~~~~~~~~~~~~~~~~~~~~~~~~~~~~~~~~~~~~~~~~~~~~
           ~~~~~\qedwhite$ 

\begin{cor} Let $R$ be a regular finitely generated ring over $\mathbb{Z}$ of dimension $d, d\geq 3, \frac{1}{d!} \in R.$ Then $E_{d+1}(R[X])$ acts transitively on $Um_{d+1}(R[X]).$ 
\end{cor}
${\pf}$ The result follows in view of $($ {\cite[Theorem 18.2]{7}}$)$ and corollary \ref{d+1}.
$~~~~~~~~~~~~~~~~~~~~~~~~~~~~~~~~~~~~~~~~~~~~~~~~~~~~~~~~~~~~~~~~~~~~~~~~~~~~~~~~~~~~~~~~~~~~~~~~~~~~~~~~~
           ~~~~~\qedwhite$ 

\begin{cor} Let $R$ be a regular finitely generated algebra over a field $K$ of dimension $d, d\geq 2, \frac{1}{d!} \in R.$ Then $E_{d+1}(R[X])$ acts transitively on $Um_{d+1}(R[X]).$ 
\end{cor}
${\pf}$ The result follows in view of $($ {\cite[Theorem 20.5]{7}}$)$ and corollary \ref{d+1}.
$~~~~~~~~~~~~~~~~~~~~~~~~~~~~~~~~~~~~~~~~~~~~~~~~~~~~~~~~~~~~~~~~~~~~~~~~~~~~~~~~~~~~~~~~~~~~~~~~~~~~~~~~~
           ~~~~~\qedwhite$ 
\vskip 0.5in

\par In $($ {\cite[Corollary 2.7]{4}}$),$ A.A. Suslin proved the following result: 
\begin{cor}$($A.A. Suslin$)$ 
Let $R$ be a noetherian ring of dimension $d$, then the canonical mapping 
$$\varphi : GL_{r}(R[X])\longrightarrow K_{1}(R[X])$$
is an epimorphism for $r\geq d+1.$
\end{cor}

As an application of corollary \ref{d+1}, we prove :
\begin{cor}$($A.A. Suslin$)$ 
Let $R$ be a regular ring of dimension $d, d\geq 2$, $\frac{1}{d!} \in R$  such that $E_{d+1}(R)$ acts transitively on $Um_{d+1}(R)$,, then the canonical mapping 
$$\varphi : GL_{r}(R[X])\longrightarrow K_{1}(R[X])$$
is an epimorphism for $r\geq d+1.$
\end{cor}
${\pf}$ The result is obvious in view of corollary \ref{d+1}.
$~~~~~~~~~~~~~~~~~~~~~~~~~~~~~~~~~~~~~~~~~~~~~~~~~~~~~~~~~~~~~~~~~~~~~~~~~~~~~~~~~~~~~~~~~~~~~~~~~~~~~~~~~
           ~~~~~\qedwhite$

\medskip
\noindent
{\bf Acknowledgement:} I thank the referee for going through
the manuscript with great care. A detailed list of suggestions by the referee improved
the exposition considerably.

\Addresses


\begin{thebibliography}{9}
\bibitem{basu} 
R. Basu;
\textit{Absence of torsion for $NK_{1}(R)$ over associative rings}, 
 Journal of Algebra and Its Applications, Vol. 10, No. 04 (2011), pp. 793--799. 


\bibitem{bloch} 
S. Bloch;
\textit{ Algebraic $K$-Theory and crystalline cohomology}, 
 Publ. Math. I.H.E.S. 47 (1977), 187--268. 


\bibitem{fasel} 
J. Fasel;
\textit{Mennicke symbols, $K$-cohomology and a Bass–Kubota theorem}, 
 Trans. Amer. Math. Soc. 367(1) (2015) 191–208.
 

\bibitem{gr} 
A. Garge, R.A. Rao;
\textit{ A nice group structure on the orbit space of unimodular rows}, 
 $K$-Theory 38 (2008), no. 2, 113--133. 
 
 \bibitem{ggr} 
A. Gupta, A. Garge, R.A. Rao; 
\textit{A nice group structure on the orbit space of unimodular rows--II}, 
 J. Algebra 407 (2014), 201--223.
 
\bibitem{mg} 
S.K. Gupta, M.P. Murthy; 
\textit{Suslin's work on linear groups over polynomial rings and Serre problem}, 
Indian Statistical Institute, New Delhi.


\bibitem{key}  S. Jose, R.A. Rao; \textit{A structure theorem for the
elementary unimodular vector group}, Trans. Amer. Math. Soc.  358
(2006), no. 7, 3097--3112.

\bibitem{orthoquot}  S. Jose, R.A. Rao; \textit{The orbit set and the orthogonal quotient},  
J. Algebra Appl. 15 (2016), no. 6, 1650101, 14 pp.

\bibitem{vdk1} 
W. van der Kallen;
\textit{A group structure on certain orbit sets of unimodular rows}, 
J. Algebra 82 (2) (1983), 363--397.

\bibitem{vdk2} 
W. van der Kallen;
\textit{A module structure on certain orbit sets of unimodular rows}, 
J. Pure Appl. Algebra 57 (1989), 657--663.

\bibitem{vdk3} 
W. van der Kallen;
\textit{From Mennicke symbols to Euler class groups}, 
Algebra, arithmetic and geometry, Part I, II (Mumbai, 2000), 341--354, Tata Inst. Fund. Res. Stud. Math., 16, 
Bombay, 2002.




\bibitem{invent}  R.A. Rao; \textit{ The Bass-Quillen conjecture in dimension three but
characteristic $\neq 2,3$ via a question of A. Suslin},
 Invent. Math. 93 (1988), no. 3, 609--618.


 
 


\bibitem{rswan}  R.A. Rao, R.G. Swan; \textit{On some actions of stably elementary matrices on alternating matrices}, 
 see excerpts on the homepage of R.G. Swan at http://www.math.uchicago.edu/~swan/.


\bibitem{st} 
J. Stienstera;
\textit{Operation in the linear $K$-theory of endomorphisms}, 
Current Trends in Algebraic Topology, Conf. Proc. Can. Math. Soc. 2 (1982).


\bibitem{sus2} 
A.A. Suslin;
\textit{On stably free modules}, 
Math. USSR-Sb. 31 (1977), no. 4, 479--491.

\bibitem{4} 
A.A. Suslin;
\textit{On the structure of special linear group over polynomial rings}, 
Math. USSR. Izv. 11 (1977), 221--238.

\bibitem{suskcohomo} 
A.A. Suslin;
\textit{$K$- theory and $K$- cohomology of certain group varieties}, 
Algebraic $K$- theory, 53--74, Adv. Soviet Math. 4 (1991), Providence, RI.



\bibitem{weibel} 
C.A. Weibel;
\textit{Mayer-Vietoris sequences and module structure on $NK_{0}$}, 
 Lecture Notes in Mathematics, Volume 854, Springer (1981), 466--498.

\bibitem{7} 
L.N. Vaserstein, A.A. Suslin;
\textit{Serre's problem on projective modules}, 
Math. USSR Izv. 10 (1976), no. 5, 937--1001.

 \bibitem{orbit} 
L.N. Vaserstein;
\textit{Operations on orbits of unimodular vectors}, 
 J. Algebra 100 (2) (1986),  456--461.

 \bibitem{V} L.N. Vaserstein; \textit{On the stabilization of the general 
linear group over a ring},  Mat. Sbornik (N.S.) 79 (121) 405--424 (Russian); 
English translation in Math. USSR Sbornik. 8 (1969), 383--400.

\end{thebibliography}
\end{document}